\documentclass[12 points]{article}
\usepackage{amssymb}
\usepackage{eucal}
\usepackage{amsmath}
\usepackage{amsthm}
\usepackage{graphicx}

\newcommand{\Log}{{\mbox{Log}}}
\newcommand{\R}{{\mathbb  R}}  \numberwithin{equation}{section} \newtheorem{thm}{\bf
Theorem}[section]
  
  \theoremstyle{remark}

\begin{document}

\title{\large\bf AVERAGING ON MANIFOLDS BY EMBEDDING ALGORITHM} \author{Petre Birtea*, Dan Com\u{a}nescu* and
C\u{a}lin-Adrian Popa**\\ {\small *Department of Mathematics, West University of Timi\c soara}\\ {\small Bd. V.
P\^ arvan,
No 4, 300223 Timi\c soara, Rom\^ania}\\ {\small birtea@math.uvt.ro, comanescu@math.uvt.ro}\\ {\small **Department of Computer Science,
"Politehnica" University of Timi\c soara}\\ {\small Bd. V. P\^ arvan, No 2, 300223 Timi\c soara, Rom\^ania}\\
{\small
calin.popa@cs.upt.ro}\\ } \date{ } \maketitle

\begin{abstract}
We will propose a new algorithm for finding critical points of cost functions defined on a differential manifold. We will lift the initial cost function to a manifold that can be embedded in a Riemannian manifold (Euclidean space) and will construct a vector field defined on the ambient space whose restriction to the embedded manifold is the gradient vector field of the lifted cost function. The advantage of this method is that it allows us to do computations in Cartesian coordinates instead of using local coordinates and covariant derivatives on the initial manifold. We will exemplify the algorithm in the case of $SO(3)$ averaging problems and will rediscover a few well known results that appear in literature.
\end{abstract}



\section{Introduction}

A problem which one frequently encounters in applications is to find an average for a finite set of sample points  $\{y_1,y_2,\ldots,y_r\}$ belonging to a manifold $N$. A method to solve this problem is to write it as an optimization problem, more precisely to construct a cost function $G_{N}:N \rightarrow \mathbb{R}$ associated with this set of sample points and then the average is defined by
$$\arg\min_{y\in N}G_N(y).$$
Of special interest are the cost functions of least square type $G_{N}(y):=\sum_{i=1}^{r}d^{2}(y,y_i)$, where $d$ is
a distance function on $N$. The average of the set of sample points $\{y_1,y_2,\ldots,y_r\}$ on the manifold $N$ is the set
defined by
$$\arg\min_{y\in N} \sum_{i=1}^{r}d^{2}(y,y_i).$$
When the function $G_{N}$ is differentiable, this is equivalent with solving the equation $dG_{N}(y)=0$
and with testing for which solutions the minimum value is attained.
In order to write this equation we need the knowledge of a local system of coordinates on the manifold $N$, a requirement that might be difficult to fulfill in many practical cases. Another problem that we may encounter with this approach is that the set of sample points $\{y_1,...,y_r\}$ might not be entirely included in the domain of a single local system of coordinates. One way to overcome these difficulties is to lift the problem on a simpler manifold $S$.

Let ${\bf P}:S\rightarrow N$ be a surjective submersion. The lifted cost function is $G_S:S\rightarrow \R$ defined by $G_S=G_N\circ{\bf P}$. The set equality ${\bf P}(\{x\in S\,|\,dG_S(x)=0\})=\{y\in N\,|\,dG_N(y)=0\}$ shows that it is sufficient to solve the equation $dG_S(x)=0$ on the simpler manifold $S$ and project these solutions on the manifold $N$.

A way to solve this new problem is to endow $S$ with a Riemannian metric $\tau$ and compute the critical points of the gradient vector field $\nabla_{\tau}G_{S}$ and verify which one of them are local/global minima. We note that if $x_0\in S$ is a critical point of the vector field $\nabla_{\tau} G_{S}$, then it is a critical point of the vector field $\nabla_{\tau'} G_{S}$, where $\tau'$ is any other Riemannian metric on $S$.
Although the manifold $S$ might have a simpler geometrical structure then the initial manifold $N$, we still need a local system of coordinates or the knowledge of covariant derivatives on the manifold $S$ in order to be able to compute the critical points of the function $G_S$, see \cite{absil-mahony-sepulchre}, \cite{absil-mahony-sepulchre-1}, \cite{edelman}, \cite{fiori-1}, \cite{moakher}, \cite{samir-absil}. To further avoid the use of Riemannian geometry of the manifold $S$ we will embed $S$ in a larger space $M$.

In the current paper we propose a solution for finding the critical points of the lifted cost function $G_S$ by constructing a vector field  ${\bf v}_0\in \mathcal{X}(M)$ on the ambient space $M$ (usually an Euclidean space), that is tangent to the submanifold $S$ and by also constructing a Riemannian metric $\tau$ on $S$ such that $\mathbf{v_0}_{|S}=\nabla _{\tau}G_S$. Consequently, the critical points of the vector field ${\bf v}_0$  (usually written in Cartesian coordinates) that belong to $S$ are critical points of the lifted cost function $G_S$ and their projection through the surjective submersion ${\bf P}$ gives the critical points of the initial cost function $G_N$. This construction will be illustrated in details for the averaging problem on the Lie group $SO(3)$ associated with four different cost functions.
Among this functions we will study the $L^p$ cost functions. We will analyze and compare the averaging problems for $L^2$ and $L^4$ cases. The $L^2$ case have been studied before in the literature. We will show that the $L^4$ case differ from the $L^2$ case in a few important aspects.

In the literature the problem of averaging on Lie groups is often solved by using the exponential map as a tool to introduce local coordinates and so lifting the problem on the tangent space, see \cite{absil-mahony-sepulchre-1}, \cite{arsigny}, \cite{fiori-1}, \cite{fiori-tanaka},  \cite{govindu}, \cite{pennec}. This method can be generalized to the context of Riemannian manifolds as there also exists an exponential map. The exponential map can be further replaced by the more general notion of retraction developed in \cite{absil-mahony-sepulchre-1}, \cite{absil-malick}, \cite{fiori}.
Averaging problems coming from real world applications have been studied in \cite{absil-mahony-sepulchre}, \cite{absil-mahony-sepulchre-1}, \cite{helmke}, \cite{samir-absil}, \cite{subbarao}, \cite{absil-vandewalle} by using the notions of covariant derivatives and the geometry of geodesics on various Riemannian manifolds.
In this paper we will propose a different algorithm in order to deal with averaging problems on general differentiable manifolds.

Other interesting cost functions are the ones that come from the Fermat-Torricelli problem. This averaging problem has been studied on various spaces, see \cite{comanescu-sever-demonstratio}, \cite{comanescu-dragomir-kikianty}, \cite{hartly}, \cite{papini-puerto}, and it will be of interest to apply the techniques developed in this paper to such problems.

\section{Averaging on a differentiable manifold}

We will solve the averaging problem presented in the Introduction by embedding the manifold $S$ as a submanifold of a Riemannian manifold $(M,g)$.
Further more, we will assume in this paper that $S$ is the preimage of a regular value for a smooth function $\mathbf{F}:=(F_1,\ldots,F_k):M\rightarrow \mathbb{R}^k$, i.e. $S=\mathbf{F}^{-1}(c_0)$, for $c_0$ a regular value of $\mathbf{F}$.

As we have already stated, our approach is to find a Riemannian metric $\tau$ on the submanifold $S$ and a vector field $\mathbf{v_0}\in\mathcal{X}(M)$ such that
\begin{equation}\label{v0-gradient-100}
\mathbf{v_0}_{|S}=\nabla_{\tau}G_S.
\end{equation}
Consequently, our initial
problem is equivalent with finding the critical points of the vector field $\mathbf{v_0}$, which is defined on
the ambient space $M$, and with choosing the ones that belong to $S$.

The vector field that we are looking for is the standard control vector field $\mathbf{v_0}$ introduced in \cite{birtea-comanescu}, which in turn  is a continuation of the study we have begun in  \cite{birtea-comanescu-siam}. Let $(M,g)$ be an $m$-dimensional Riemannian manifold and $F_1,\ldots,F_k,G:M\rightarrow\mathbb{R}$ be $k+1$ smooth functions. The standard control vector field is defined by
\begin{equation}\label{v0}
    \mathbf{v_0}=\sum_{i=1}^k(-1)^{i+k+1}\det \Sigma_{(F_1,\ldots ,\widehat{F_i},\ldots ,F_k,G)}^{(F_1,\ldots
    ,F_k)}\nabla
    F_i+\det\Sigma_{(F_1,\ldots ,F_k)}^{(F_1,\ldots ,F_k)}\nabla G,
\end{equation}
where  $\widehat{\cdot}$ represent the missing term and $\Sigma_{(g_1,...,g_s)}^{(f_1,...,f_r)}$ is the $r\times s$ Gramian matrix
\begin{equation}\label{sigma}
\Sigma_{(g_1,...,g_s)}^{(f_1,...,f_r)}=\left(%
\begin{array}{cccc}
  <\nabla g_1,\nabla f_{1}> & ... & <\nabla g_s,\nabla f_{1}> \\
  ... & ... & ... \\

  <\nabla g_1,\nabla f_r> & ... & <\nabla g_s,\nabla f_r> \\
\end{array}%
\right),
\end{equation}
generated by the smooth functions $f_1,...,f_r,g_1,...,g_s:M\rightarrow \mathbb{R}$.

The vector field $\mathbf{v_0}$ conserves the regular leaves
of the map $\mathbf{F}:=(F_1,\ldots ,F_k):M\rightarrow \mathbb{R}^k$ and dissipates the function $G$, more
precisely the derivation of the function $G$ along the vector field ${\bf v}_0$ is given by the Lie derivative $L_{\mathbf{v_0}}G = \det\Sigma_{(F_1,\ldots ,F_k,G)}^{(F_1,\ldots ,F_k,G)}\geq 0$.

In what follows, we will describe the geometry of the standard control vector field.
Let $\Omega^1(M)$ be the real vector space of the differential one forms on the manifold $M$. Let $\mathbf{T}:\Omega^1(M)\times\Omega^1(M)\rightarrow \mathbb{R}$ be the degenerate  symmetric contravariant
2-tensor given by
\begin{equation}\label{T}
\mathbf{T}:=\sum_{i,j=1}^k(-1)^{i+j+1}\det\Sigma_{(F_1,\ldots ,\hat{F_i},\ldots ,F_k)}^{(F_1,\ldots
,\hat{F_j},\ldots ,F_k)}\nabla F_i\otimes\nabla F_j+\det\Sigma_{(F_1,\ldots ,F_k)}^{(F_1,\ldots ,F_k)} g^{-1},
\end{equation}
where $g^{-1}$ is the cometric 2-tensor $g^{-1}(x)=g^{pq}(x)\frac{\partial}{\partial x^p}\otimes \frac{\partial}{\partial x^q}$ constructed from the metric tensor $g$ and the contravariant 2-tensor
$\nabla F_i\otimes\nabla F_j :\Omega^1(M)\times\Omega^1(M)\rightarrow \mathbb{R}$  is defined by the formula
$$\nabla F_i\otimes\nabla F_j\,\,(\alpha,\beta):=\alpha (\nabla F_i)\beta(\nabla F_j).$$

In Riemannian geometry the gradient vector field of a smooth function $G$ can be defined by the formula $\nabla G={\bf i}_{dG}g^{-1}$, where ${\bf i}$ denotes the interior product defined by ${\bf i}_{dG}g^{-1}(\alpha):=g^{-1}(dG,\alpha),\,\,\alpha\in \Omega^1(M)$.
The following result shows that $\mathbf{v_0}$ looks like the gradient vector field of the function $G$, with
respect to the degenerate symmetric contravariant 2-tensor $\mathbf{T}$.

\begin{thm}\textit{{\cite{birtea-comanescu}}}\label{v0-cu-i}
On the manifold $(M,g)$, the standard control vector field $\mathbf{v_0}$ is given by the following formula
$$\mathbf{v_0}=\mathbf{i}_{dG}\mathbf{T}.$$
\end{thm}
From the above theorem, computing the critical points of $\mathbf{v_0}$ (which is written in local coordinates on $M$) that belong also to $S$ is
equivalent with solving the following system of $m+k$ equations on $M$:
\begin{equation}\label{sistem}
\left\{
\begin{array}{ll}
\left[%
\begin{array}{ccccccc}
   \mathbf{T}^{11}(x) & \ldots & \mathbf{T}^{1m}(x) \\
   \vdots & \cdots & \vdots \\
    \mathbf{T}^{m1}(x)& \ldots & \mathbf{T}^{mm}(x) \\
\end{array}%
\right]\left[%
\begin{array}{c}
  \displaystyle\frac{\partial G}{\partial x^{1}}(x)\\ \vdots \\ \displaystyle\frac{\partial G}{\partial
  x^{m}}(x)
\end{array}%
\right]={\bf 0}\\
\mathbf{F}(x)=c
\end{array}
\right.,
\end{equation}
where $\displaystyle\mathbf{T}=\mathbf{T}^{pq}\frac{\partial}{\partial x^{p}}\otimes\frac{\partial}{\partial
x^{q}}$ and $\displaystyle dG=\frac{\partial G}{\partial x^{i}}dx^i$.
In general, a system of $m+k$ equations with $m$ unknowns might not have any solutions. But the above system
that has $m$ unknowns, namely the coordinates of $x$ in $M$, has exactly $m$ functional independent equations
due to the fact that the rank of the tensor $\mathbf{T}$ is $m-k$ for every regular point in the open set in
$M$ of regular points of $\mathbf{F}$ (the rank of the tensor ${\bf T}$ is $m-k$ at every point $x\in M$ where $\nabla F_1(x),...,\nabla F_k(x)$ are linear independent, see Section 4 of \cite{birtea-comanescu}).

The standard control vector field ${\bf v_0}$ is tangent to every regular leaf and next we will recall its geometry when restricted to a regular leaf.
Let $i_c:L_c\rightarrow M$  be the canonical inclusion of the regular leaf $L_c$ into the manifold $M$. The
contravariant symmetric 2-tensor $\mathbf{T}$ is non-degenerate when restricted to a regular leaf $L_c$ and
consequently, the restriction can be inverted and it thus generates the Riemannian metric, see
\cite{birtea-comanescu}, $\tau_c:\mathcal{X}(L_c)\times \mathcal{X}(L_c)\rightarrow \mathcal{C}^{\infty}(L_c)$
defined by
$$\tau_c(X^c,Y^c):=\mathbf{T}^{-1}((i_{c})_{*}X^c,(i_{c})_{*}Y^c),$$
where $(i_{c})_{*}:\mathcal{X}(L_c)\rightarrow \mathcal{X}(M)$ is the push-forward application.

\begin{thm}\cite{birtea-comanescu}\label{caracterizare}
On a regular leaf $L_c$ we have the following characterizations.
\begin{itemize}
\item [(i)] $\tau_c=\frac{1}{\det\Sigma_{(F_1,\ldots ,F_k)}^{(F_1,\ldots ,F_k)}\circ i_c}(i_c)^*g$, where $(i_c)^*$ is the pull-back operator;
\item [(ii)] $\mathbf{v_0}_{|L_c}=(i_{c})_{*}\nabla_{\tau_c}(G\circ i_c)$.
\end{itemize}
\end{thm}

The above theorem shows that the tensor ${\bf T}$ induces a Riemannian metric on every regular leaf $L_c$ which is the first fundamental form of the submanifold $L_c\subset (M,g)$ multiplied with a positive function. The vector field ${\bf v}_0$ is tangent to the submanifold $L_c$ and the restriction is equal with the gradient vector field of the restricted function $G_{|L_c}$ with respect to the metric $\tau_c$.
As $S$ is such a regular leaf, the Riemannian metric on $S$ is $\tau=\tau_{c_0}$ and the above theorem shows that the equality \eqref{v0-gradient-100} holds, where $G_S=G_{|S}$.
\medskip

In conclusion, in order to solve our initial averaging problem on the manifold $N$ associated to the cost function $G_N$, by using the standard control vector field, we apply the following {\bf embedding algorithm}:
\begin{itemize}
\item [{\bf (i)}] Choose the geometrical setting of a manifold $S$, a surjective submersion ${\bf P}:S\rightarrow N$ and construct the lifted cost function $G_S=G_N\circ {\bf P}$.

\item [{\bf (ii)}] Find a Riemannian ambient space $(M,g)$ and a differentiable function ${\bf F}:M\rightarrow \R^k$ such that $S={\bf F}^{-1}(c_0)$, where $c_0$ is a regular value of ${\bf F}$. Construct the contravariant symmetric 2-tensor $\mathbf{T}$.

\item [{\bf (iii)}] Find a prolongation function $G:M\rightarrow \R$ such that $G_{|S}=G_S$ and construct the standard control vector field ${\bf v_0}$ with the initial data ${\bf F},G$.

\item [{\bf (iv)}] Solve the system $\mathbf{v_{0}}_{|S}(x)={\bf 0}$ (i.e. the system \eqref{sistem}) which is equivalent with finding the critical points of the lifted cost function $G_S$.
Project through ${\bf P}$ these solutions and thus obtain the critical points of the initial cost function $G_N$.

\item [{\bf (v)}] Find critical points that are local/global minima for the cost function $G_N$.
\end{itemize}

Examples that fit the geometry underlying this embedding algorithm are given by the rich cases of quotient manifolds, $N=S/K$ with $K$ a group that acts on $S$.
In the case when $N$ can be written as a preimage of a regular value of a smooth map ${\bf F}$, then {\bf (i)} is not a necessary step of the algorithm ($S=N$ and ${\bf P}$ is the identity).

In the case when ${\bf P}_*{\bf v_0}_{|S}$ is a vector field on the manifold $N$, then we have the set inclusion ${\bf P}(\{x\in S\,|\,{\bf v_0}_{|S}(x)={\bf 0}\})\subset \{y\in N\,|\,{\bf P}_*{\bf v_0}_{|S}(y)={\bf 0}\}$. If ${\bf P}$ is also a local diffeomorphism we will obtain equality between the above two sets and consequently, the critical points of the cost function $G_N$ are characterized by the equation ${\bf P}_*{\bf v_0}_{|S}(y)={\bf 0}.$ In this particular setting the step {\bf (iv)} can be replaced by the following:

\begin{itemize}
\item [{\bf (iv')}] Suppose that  ${\bf P}_*{\bf v_0}_{|S}\in \mathcal{X}(N)$ and ${\bf P}$ is a local diffeomorphism, then solve the equation
${\bf P}_*{\bf v_0}_{|S}(y)={\bf 0}$ (the solutions are critical points for the cost function $G_N$).
\end{itemize}

For solving point {\bf (v)} one can study the second order derivatives of the cost function $G_N$ (see \cite{edelman}, \cite{helmke}, \cite{moakher}) or study the convexity properties of the cost functions as in \cite{hartly}.

\section{Averaging on $SO(3)$}
A problem frequently encountered in practice is finding a rotation that is in some sense an average of a finite set of rotations. The group of special rotations is not a Euclidean space, but a differential manifold. Consequently, an averaging problem on such a space has to be considered in the realm of differential geometry.

The special orthogonal group $SO(3)$ can be identified by a double covering map with the sphere $S^{3}$ in $\mathbb{R}^{4}$, see \cite{altman}, \cite{shuster}. Our manifold $S$ will be $S^{3}$ and the ambient space $(M,g)$ will be the Euclidean space $\mathbb{R}^{4}$. We will consider different distance functions on $SO(3)$ that generate different cost functions and will solve the associated averaging problem, as we have described in the previous section, by using the embedding algorithm .

From a historical perspective, we will use the quaternion notation on $\R^4$, although the quaternion structure will not be specifically used in our computations.  The unit quaternions $\mathbf{q}=(q^0,q^1,q^2,q^3)\in S^3\subset \mathbb{R}^{4}$ and $-\mathbf{q}\in S^3\subset\mathbb{R}^{4}$ correspond to the following rotation in $SO(3)$:
\begin{equation*}
\mathbf{R}^{\mathbf{q}}=\left(
\begin{array}{ccc}
(q^0)^{2}+(q^{1})^{2}-(q^{2})^{2}-(q^{3})^{2} & 2(q^{1}q^{2}-q^{0}q^{3}) & 2(q^{1}q^{3}+q^{0}q^{2}) \\
2(q^{1}q^{2}+q^{0}q^{3})  & (q^0)^{2}-(q^{1})^{2}+(q^{2})^{2}-(q^{3})^{2}& 2(q^{2}q^{3}-q^{0}q^{1}) \\
2(q^{1}q^{3}-q^{0}q^{2}) & 2(q^{2}q^{3}+q^{0}q^{1}) & (q^0)^{2}-(q^{1})^{2}-(q^{2})^{2}+(q^{3})^{2}
\end{array}
\right)
\end{equation*}

This gives rise to the smooth double covering map $\mathbf{P}:S^{3}\rightarrow SO(3)$,
$\mathbf{P}(\mathbf{q})=\mathbf{R}^{\mathbf{q}}$. The covering map is a local diffeomorphism and consequently,
instead of working with the distance function (cost function) on $SO(3)$, we will work with cost functions on the unit
sphere $S^{3}\subset \mathbb{R}^{4}$ as it has been explained in the previous section. In this case, for the step {\bf (ii)} of the embedding algorithm, we have $S^{3}=F^{-1}(1)$, with
$F:\mathbb{R}^{4}\rightarrow\mathbb{R}$,
$F(\mathbf{q})=(q^{0})^{2}+(q^{1})^{2}+(q^{2})^{2}+(q^{3})^{2}$. The tensor ${\bf T}$ is given by the formula, see \cite{birtea-comanescu} eq. (4.2), ${\bf T}=-\nabla F\otimes \nabla F+||\nabla F||^2 g^{-1}$, where $g$ is the Euclidean metric  on $\mathbb{R}^4$.
Writing this explicitly, the matrix associated to the symmetric contravariant 2-tensor
$\mathbf{T}$ is given by:
$$[\mathbf{T}(\mathbf{q})]=4
\left(
  \begin{array}{cccc}
    (q^{1})^{2}+(q^{2})^{2}+(q^{3})^{2} & -q^{0}q^{1} & -q^{0}q^{2} & -q^{0}q^{3} \\
    -q^{1}q^{0} & (q^{0})^{2}+(q^{2})^{2}+(q^{3})^{2} & -q^{1}q^{2} & -q^{1}q^{3} \\
    -q^{2}q^{0} & -q^{2}q^{1} & (q^{0})^{2}+(q^{1})^{2}+(q^{3})^{2} & -q^{2}q^{3} \\
    -q^{3}q^{0} & -q^{3}q^{1} & -q^{3}q^{2} & (q^{0})^{2}+(q^{1})^{2}+(q^{2})^{2}
  \end{array}
\right).$$

In order to compute the standard control vector field ${\bf v}_0$ according to the formula given in Theorem \ref{v0-cu-i}, we need the following elementary computation:
if $\omega({\bf q})=\omega_0({\bf q})dq^0+\omega_1({\bf q})dq^1+\omega_2({\bf q})dq^2+\omega_3({\bf q})dq^3$ is a one form on $\R^4$, then we have the vector field
\begin{equation}\label{formula-calcul}
{\bf i}_{\omega}{\bf T}({\bf q})=4(\langle {\bf q},{\bf q} \rangle \overline{\omega}({\bf q})-\langle {\bf q},\overline{\omega}({\bf q}) \rangle {\bf q}),
\end{equation}
where we have made the notation $\overline{\omega}({\bf q}):=(\omega_0({\bf q}),\omega_1({\bf q}),\omega_2({\bf q}),\omega_3({\bf q}))$.

We will exemplify our construction by using various distance functions used in literature for measuring
distances between Euclidean transformations. An extensive list is presented in \cite{huynh}, where their equivalence and functional dependence is also studied.
\medskip

{\bf I.}
On the Lie group $SO(3)$ we will consider the distance function $d_1:SO(3)\times
SO(3)\rightarrow\mathbb{R}_{+}$, $d_1(\mathbf{R}_1,\mathbf{R}_2)=||\mathbf{R}_1-\mathbf{R}_2||_{F}$, where $||\cdot||_F$ is the Frobenius norm. The associated cost function is $G_{1_{SO(3)}}:SO(3)\rightarrow \R$,
$$G_{1_{SO(3)}}({\bf R})=\sum_{i=1}^r||{\bf R}-{\bf R}_i||_F^2,$$
where ${\bf R}_1,...,{\bf R}_r$ are the sample rotation matrices,  see \cite{moakher}, \cite{sharf-wolf-rubin}.

According to the step {\bf (i)} of the embedding algorithm,
we will lift the problem of finding the critical points of the cost function $G_{1_{SO(3)}}$ to the problem of finding the critical points of the lifted cost function $G_{{1}_{S^{3}}}:S^3\rightarrow \R$, $G_{{1}_{S^{3}}}:=G_{1_{SO(3)}}\circ {\bf P}$.

In order to compute the lifted cost function $G_{{1}_{S^{3}}}$, by using the surjectivity of the covering map ${\bf P}$, for each sample matrix ${\bf R}_i$ we  will choose a sample quaternion ${\bf q}_i\in S^3$.
We have the following computation:
\begin{align*}
G_{{1}_{S^{3}}}({\bf q}) & = \sum_{i=1}^r||\mathbf{R^{q}}-\mathbf{R^{q_{i}}}||_{F}^2 =  \sum_{i=1}^r \text{tr}(({\bf R}^{\bf q}-{\bf R}^{\bf q_i})({\bf R}^{\bf q}-{\bf R}^{\bf q_i})^T) \\
& = 2\sum_{i=1}^{r}(3-\text{tr}({{\bf R}^{\bf q}}^T{\bf R}^{\bf q_i})) \\
& =8\sum_{i=1}^{r}
(1-\langle\mathbf{q},\mathbf{q}_{i}\rangle^{2}),
\end{align*}
where we have used the equality
\begin{equation}\label{modul}
\langle\mathbf{q},\mathbf{q}_{i}\rangle ^2=\frac{1}{4}(\text{tr}({{\bf R}^{\bf q}}^T{\bf R}^{{\bf q}_i})+1).
\end{equation}

For implementing the step {\bf (iii)} we need to construct the prolongation function
$G_{1}:\mathbb{R}^{4}\rightarrow\mathbb{R}$, $G_{1}(\mathbf{q})=\sum_{i=1}^{r}
8(1-\langle\mathbf{q},\mathbf{q}_{i}\rangle^{2})$. Note that $G_{{1}_{S^{3}}}$ and $G_1$ are even functions so they do not depend on the choice of the sample quaternions $\mathbf{q}_i$ or $-\mathbf{q}_i$ which represent the same sample rotation ${\bf R}_i$. This subtler problem has also been addressed in a different way in \cite{markley}.
Computing the differential 1-form $dG_1$ we obtain:
\begin{align*}
\overline{d G_{1}}(\mathbf{q}) & =\left(\frac{\partial G_1}{\partial q^{0}}(\mathbf{q}),...,\frac{\partial G_1}{\partial
  q^{3}}(\mathbf{q})
\right)=\left(
  -16\sum_{i=1}^{r}q_i^0<\mathbf{q},\mathbf{q}_i>,
...,
  -16\sum_{i=1}^{r}q_i^3<\mathbf{q},\mathbf{q}_i>
\right) \\
& =-16\sum_{i=1}^{r}<\mathbf{q},\mathbf{q}_i>\mathbf{q}_i.
\end{align*}
Using \eqref{formula-calcul}, where the 1-form $\omega$ is $dG_1$,
the system of equations \eqref{sistem} becomes:
\begin{equation}\label{sistem-v0-I}
\left\{\begin{array}{ll}
\displaystyle\sum_{i=1}^{r}\langle\mathbf{q},\mathbf{q}_{i}\rangle (\langle\mathbf{q},\mathbf{q}\rangle\mathbf{q}_{i}-\langle\mathbf{q},\mathbf{q}_{i}\rangle\mathbf{q})
={\bf 0}
\\<\mathbf{q},\mathbf{q}>=1
\end{array}
\right.,
\end{equation}
and this represents the equation $\mathbf{v_{0}}_{|S^{3}}(\mathbf{q})={\bf 0}$.
The above system has four equations of third degree in the four unknowns $\mathbf{q}=(q^0,q^1,q^2,q^3)$ corresponding to the standard control vector field $\mathbf{v_0}={\bf 0}$ plus the constraint equation that describes $S^3$.

The four equations corresponding to $\mathbf{v_0}={\bf 0}$ are not functionally independent
because, as we have stated in Section 2, the tensor $\mathbf{T}$ is degenerate. Nevertheless, since
$\mathbf{v_{0}}(\mathbf{q})\in T_{\mathbf{q}}S^{3}$, for any $\mathbf{q}\in S^3$, the system
\eqref{sistem-v0-I} can be described only by three equations plus the constraint equation.

In order to write $\mathbf{v_{0}}_{|S^{3}}(\mathbf{q})={\bf 0}$ as a system of three independent equations, we will push forward the vector field $\mathbf{v_{0}}_{|S^{3}}$ through the covering map
$\mathbf{P}$, which will also  place us in the hypotheses of the step {\bf (iv')} of the embedding algorithm. By definition of the push-forward operator, we have that $$\mathbf{P_*}\mathbf{v_{0}}_{|S^{3}}(\mathbf{R}^{\mathbf{q}})=
D\mathbf{P}({\mathbf{q}})\cdot \mathbf{v_{0}}_{|S^{3}}(\mathbf{q})\in T_{\mathbf{R}^{\mathbf{q}}}SO(3).$$
As $\mathbf{P}$ is not an injective map, in order to obtain a vector field on the target space $SO(3)$, the equality
$D\mathbf{P}({\mathbf{q}})\cdot \mathbf{v_{0}}_{|S^{3}}(\mathbf{q})=D\mathbf{P}(-{\mathbf{q}})\cdot \mathbf{v_{0}}_{|S^{3}}(-\mathbf{q})$
must be satisfied. We will show in what follows that this is the case.

We need to compute the tangent map of the
covering map $\mathbf{P}:S^3\rightarrow SO(3)$. This map can be written as a restriction of the map $\widetilde{\mathbf{P}}:\R^4\rightarrow \R^9$ defined by $\widetilde{\mathbf{P}}(q^0,q^1,q^2,q^3)=((q^0)^{2}+(q^{1})^{2}-(q^{2})^{2}-(q^{3})^{2}, 2(q^{1}q^{2}-q^{0}q^{3}), 2(q^{1}q^{3}+q^{0}q^{2}),
2(q^{1}q^{2}+q^{0}q^{3}), (q^0)^{2}-(q^{1})^{2}+(q^{2})^{2}-(q^{3})^{2}, 2(q^{2}q^{3}-q^{0}q^{1}),
2(q^{1}q^{3}-q^{0}q^{2}), 2(q^{2}q^{3}+q^{0}q^{1}), (q^0)^{2}-(q^{1})^{2}-(q^{2})^{2}+(q^{3})^{2})$, where
a rotation matrix $\mathbf{R}^{\mathbf{q}}$ has been identified with a point in $\mathbb{R}^{9}$ (the first line is identified with the first three components and so on).

The matrix corresponding to the linear map $D\widetilde{\mathbf{P}}({\bf q}):\R^4\rightarrow
\R^9$ is given by Jacobian matrix
$$2
\left[\begin{array}{cccc}
  q^0 &  q^1 & -q^2 & -q^3 \\
 -q^3 &  q^2 &  q^1 & -q^0 \\
  q^2 &  q^3 &  q^0 &  q^1 \\
  q^3 &  q^2 &  q^1 &  q^0 \\
  q^0 & -q^1 &  q^2 & -q^3 \\
 -q^1 & -q^0 &  q^3 &  q^2 \\
 -q^2 &  q^3 & -q^0 &  q^1 \\
  q^1 &  q^0 &  q^3 &  q^2 \\
  q^0 & -q^1 & -q^2 &  q^3
\end{array}\right].$$
By direct computation, where we identify a vector in $\R^9$ with a tangent vector $\mathbf{R}^{\mathbf{q}}\Delta\in T_{\mathbf{R}^{\mathbf{q}}}SO(3)=\{\mathbf{R}^{\mathbf{q}}\Delta\,|\,\Delta\,\text{is a}\,3\times 3\, \text{skew symmetric matrix}\}$, for $\mathbf{q}\in S^3$ we have that
\begin{align*}
D\mathbf{P}({\mathbf{q}})\cdot {{\bf v}_{0}}_{|S^3}({\bf q}) & =D\widetilde{{\bf P}}({\bf q})\cdot {{\bf v}_{0}}_{|S^3}({\bf q}) \\
& = \sum_{i=1}^{r} \langle\mathbf{q},\mathbf{q}_{i}\rangle D\widetilde{{\bf P}}({\bf q})\cdot (\langle\mathbf{q},\mathbf{q}\rangle{\bf q}_i-\langle\mathbf{q},\mathbf{q}_i\rangle{\bf q}) \\
& \cong \sum_{i=1}^{r} \langle\mathbf{q},\mathbf{q}_{i}\rangle {\bf R}^{\bf q}\Delta_i({\bf q}) \\
& = {\bf R}^{\bf q} \sum_{i=1}^{r}\langle\mathbf{q},\mathbf{q}_{i}\rangle \Delta_i({\bf q}),
\end{align*}
where
$$\Delta_{i}(\mathbf{q})=\left(\begin{array}{ccc}
  0 &  -q^{0}q_{i}^{3}+q^{1}q_{i}^{2}-q^{2}q_{i}^{1}+q^{3}q_{i}^{0} &
  q^{0}q_{i}^{2}+q^{1}q_{i}^{3}-q^{2}q_{i}^{0}-q^{3}q_{i}^{1}\\
 q^{0}q_{i}^{3}-q^{1}q_{i}^{2}+q^{2}q_{i}^{1}-q^{3}q_{i}^{0} &  0 &
 -q^{0}q_{i}^{1}+q^{1}q_{i}^{0}+q^{2}q_{i}^{3}-q^{3}q_{i}^{2}  \\
   -q^{0}q_{i}^{2}-q^{1}q_{i}^{3}+q^{2}q_{i}^{0}+q^{3}q_{i}^{1} &
   q^{0}q_{i}^{1}-q^{1}q_{i}^{0}-q^{2}q_{i}^{3}+q^{3}q_{i}^{2} &  0
\end{array}\right).$$
By identifying the vector $D\widetilde{\mathbf{P}}({\bf q})\cdot (\langle\mathbf{q},\mathbf{q}\rangle{\bf q}_i-\langle\mathbf{q},\mathbf{q}_i\rangle{\bf q})\in \R^9$ with a $3\times 3$ matrix,
the skew-symmetric matrix is given by the formula
$$\Delta_{i}(\mathbf{q})={{\bf R}^{\bf q}}^TD\widetilde{\mathbf{P}}({\bf q})\cdot (\langle\mathbf{q},\mathbf{q}\rangle{\bf q}_i-\langle\mathbf{q},\mathbf{q}_i\rangle{\bf q}).$$
By noticing that $\Delta_i(-{\bf q})=-\Delta_i({\bf q})$, we obtain the equality $D\mathbf{P}({\bf q})\cdot \mathbf{v_{0}}_{|S^{3}}(\mathbf{q})=D\mathbf{P}(-{\bf q})\cdot \mathbf{v_{0}}_{|S^{3}}(-\mathbf{q})$
and consequently, $\mathbf{P_*}\mathbf{v_{0}}_{|S^{3}}$ is a vector field on the manifold $SO(3)$.

As a result we obtain that the critical points of the cost function $G_{{1}_{S^{3}}}$ are the solutions of the following system (which is equivalent with \eqref{sistem-v0-I}):
\begin{equation}\label{terta}
\left\{
\begin{array}{ll}
    \displaystyle\sum_{i=1}^{r}\langle\mathbf{q},\mathbf{q}_{i}\rangle\Delta_{i}(\mathbf{q})= {\bf 0} \\
    \langle\mathbf{q},\mathbf{q}\rangle=1
\end{array}
\right..
\end{equation}
The advantage of the above system, compared with the system \eqref{sistem-v0-I}, is that we have only three equations of second degree instead of four equations of third degree plus the constraint equation.

If we transform the above system into rotations, by direct computation, we obtain:
\begin{equation}\label{relatie-delta}
\displaystyle\langle\mathbf{q},\mathbf{q}_{i}\rangle\Delta_{i}(\mathbf{q})=\frac{1}{4}({\mathbf{R}^{\mathbf{q}_{i}}}^{T}\mathbf{R}^{\mathbf{q}}-
{\mathbf{R}^{\mathbf{q}}}^{T}\mathbf{R}^{\mathbf{q}_{i}})
\end{equation}
 and the above system becomes:
\begin{equation*}
\sum_{i=1}^{r}({\mathbf{R}^{\mathbf{q}_{i}}}^{T}\mathbf{R}^{\mathbf{q}}-
{\mathbf{R}^{\mathbf{q}}}^{T}\mathbf{R}^{\mathbf{q}_{i}})={\bf 0}.
\end{equation*}
This is equivalent with the following equation on $SO(3)$:
\begin{equation}\label{media-aritmetica}
\overline{\bf R}^T \mathbf{R}-{\mathbf{R}}^T\overline{\bf R}={\bf 0},
\end{equation}
where $\overline{\bf R}=\frac{1}{r}\sum_{i=1}^{r}\mathbf{R}^{\mathbf{q}_{i}}=\frac{1}{r}\sum_{i=1}^{r}\mathbf{R}_i$. Solving this equation corresponds to step {\bf (iv')} of the embedding algorithm and
the solutions are the critical points of the cost function  $G_{{1}_{SO(3)}}$. Also, this equation is the same as the characterization for the Euclidean mean obtained in \cite{moakher}.
Thus, we have obtained that searching for critical points of the cost function $G_{{1}_{SO(3)}}$ can be performed in two ways. More precisely, we can solve \eqref{media-aritmetica} or we can transform the initial problem into quaternions and solve the system \eqref{terta} and then transform the solutions into rotations.

\medskip

{\bf II.}
Next, we will consider the geodesic distance on the Lie group $SO(3)$ which is defined by the angle of two rotations, see \cite{altman}, \cite{arsigny}, \cite{gramkow}, \cite{moakher}, \cite{pennec}, \cite{sharf-wolf-rubin}. For two rotations ${\bf R}_1, {\bf R}_2\in SO(3)$, $d_2({\bf R}_1,{\bf R}_2)=||\Log ({\bf R}_1^T{\bf R}_2)||_F=\sqrt{2}|\theta |$. where $\theta\in (-\pi,\pi)$ is the angle between the rotations  ${\bf R}_1$ and ${\bf R}_2$.
The associated cost function is $G_{2_{SO(3)}}:SO(3)\backslash \mathcal{A}\rightarrow \R$,
$$G_{2_{SO(3)}}({\bf R})=\sum_{i=1}^r ||\Log ({\bf R}_i^T{\bf R})||_F^2,$$
where ${\bf R}_1,...,{\bf R}_r$ are the sample rotation matrices and the closed set $\mathcal{A}$ is defined by  $\mathcal{A}=\bigcup_{i=1}^r\mathcal{A}_i$ with
$\mathcal{A}_i=\{{\bf R}\in SO(3)\,|\,\text{tr}({\bf R}_i^T{\bf R})=-1\}.$

For $\Log(\cdot)$ to be well defined we have to remove the rotations included in the closed set $\mathcal{A}$. When we lift to quaternions we use the set equality $\mathcal{A}_i={\bf P}(\Pi_i)$, where  $\Pi_i:=\{{\bf q}\in \R^4\,|\,\langle {\bf q},{\bf q}_i\rangle =0\}$. Consequently, from $S^3$ we have to subtract the set $\bigcup_{i=1}^r\Pi_i$.  As before, we will compute the lifted cost function $G_{{2}_{S^{3}}}: S^3\backslash \bigcup_{i=1}^r\Pi_i \rightarrow \R$ as stated in step {\bf (i)} of the embedding algorithm. More precisely, we have\footnote{$|\theta_i|=\arccos\left (\frac{\text{tr} ( \mathbf{R^{q_i}}^T\mathbf{R^{q}})-1}{2}\right ),\,\,\,\text{if}\,\,x\in [-1,1]\,\,\text{then}\,\,\arccos(2x^2-1)=2\arccos(|x|)  $.}
\begin{align*}
G_{{2}_{S^{3}}}({\bf q}) & = \sum_{i=1}^r||\Log(\mathbf{R^{q_i}}^T\mathbf{R^{q}})||_{F}^2 =  2\sum_{i=1}^r \theta_i^2 \\
& = 2\sum_{i=1}^{r}\arccos^2\left ( \frac{\text{tr}(\mathbf{R^{q_i}}^T\mathbf{R^{q}})-1}{2}\right ) \\
& \stackrel{\eqref{modul}}{=}2\sum_{i=1}^{r}\arccos^2(2\langle {\bf q}_i,{\bf q}\rangle^2-1) \\
& =2\sum_{i=1}^{r}\arccos^2(|\langle {\bf q}_i,{\bf q}\rangle|) .
\end{align*}

 The lifted cost function $G_{{2}_{S^{3}}}$ does not depend on the choice of sample quaternions ${\bf q}_i$ or $-{\bf q}_i$ that represent the same sample rotation ${\bf R}_i$.

For step {\bf (iii)} of the embedding algorithm, the prolongation function is given by $G_{2}:\mathbb{R}^{4}\backslash\{{\bf 0}\}\rightarrow\mathbb{R}$,
$G_{2}(\mathbf{q})=2\displaystyle\sum_{i=1}^{r}
\arccos^{2}\left(\frac{|\langle\mathbf{q},\mathbf{q}_{i}\rangle|
}{||\mathbf{q}||\cdot ||\mathbf{q}_{i}||}\right)$.
The coefficients of the differential 1-form $dG_2$ are given by:
\begin{align*}
\overline{d G_{2}}(\mathbf{q}) & =\left(\frac{\partial G_2}{\partial q^{0}}(\mathbf{q}),...,\frac{\partial G_2}{\partial
  q^{3}}(\mathbf{q})
\right) \\
& =\left(...,
   -4\sum_{i=1}^{r}\frac{\text{sgn}(\langle\mathbf{q},\mathbf{q}_{i}\rangle)||\mathbf{q}||^{2}q_{i}^{j}-|\langle\mathbf{q},\mathbf{q}_{i}\rangle|q^{j}}{||\mathbf{q}||^{2}
  \sqrt{||\mathbf{q}||^{2}||\mathbf{q}_{i}||^{2}-\langle\mathbf{q},\mathbf{q}_{i}\rangle^{2}}}
  \arccos\left(\frac{|\langle\mathbf{q},\mathbf{q}_{i}\rangle|}{||\mathbf{q}||\cdot ||\mathbf{q}_{i}||}\right),
... \right) \\
& = -4\sum_{i=1}^{r} \arccos\left(\frac{|\langle\mathbf{q},\mathbf{q}_{i}\rangle|}{||\mathbf{q}||\cdot ||\mathbf{q}_{i}||}\right)\frac{\text{sgn}(\langle\mathbf{q},\mathbf{q}_{i}\rangle)}{||\mathbf{q}||^{2}
  \sqrt{||\mathbf{q}||^{2}||\mathbf{q}_{i}||^{2}-\langle\mathbf{q},\mathbf{q}_{i}\rangle^{2}}}(\langle\mathbf{q},\mathbf{q}\rangle\mathbf{q}_{i}-\langle\mathbf{q},\mathbf{q}_{i}\rangle\mathbf{q}).
\end{align*}

In order to solve the differentiability problem we will eliminate from the domain of definition of the function $G_2$ the hyperplanes $\Pi_i=\{{\bf q}\in \R^4\,|\,\langle {\bf q},{\bf q}_i\rangle =0\}$ and the lines $d_i$ that go through points $\mathbf{0}$ and ${\mathbf q}_i$.

In this case, the equation $\mathbf{v_0}_{|S^3}(\mathbf{q})={\bf 0}$ is equivalent with
$$\left\{\begin{array}{ll}
\displaystyle\sum_{i=1}^{r}\frac{\text{sgn}(\langle {\bf q},{\bf q}_i\rangle)\arccos(|\langle\mathbf{q},\mathbf{q}_{i}\rangle|)}{\sqrt{1-\langle\mathbf{q},\mathbf{q}_{i}\rangle^{2}}}
(\langle\mathbf{q},\mathbf{q}\rangle\mathbf{q}_{i}-\langle\mathbf{q},\mathbf{q}_{i}\rangle\mathbf{q})
={\bf 0}
\\ \displaystyle {{\bf q}\in S^3\backslash \bigcup_{i=1}^r\Pi_i} \\
{\bf q}\neq \pm{\bf q}_i
\end{array}
\right.$$

By using the same arguments as before, the equation
$\mathbf{v_{0}}_{|S^3}(\mathbf{q})={\bf 0}$ is equivalent with $\mathbf{P_*}\mathbf{v_{0}}_{|S^3}={\bf 0}$ which has the following form
\begin{equation}\label{secunda}
\left\{
\begin{array}{ll}
    \displaystyle\sum_{i=1}^{r}\frac{\text{sgn}(\langle {\bf q},{\bf q}_i\rangle)\arccos(|\langle\mathbf{q},\mathbf{q}_{i}\rangle|)}{\sqrt{1-\langle\mathbf{q},\mathbf{q}_{i}\rangle^{2}}}\Delta_{i}(\mathbf{q})
    =
   {\bf  0} \\
  \displaystyle{{\bf q}\in S^3\backslash \bigcup_{i=1}^r\Pi_i} \\
{\bf q}\neq \pm{\bf q}_i
\end{array}
\right..
\end{equation}

The above expression also shows that $\mathbf{P_*}\mathbf{v_{0}}_{|S^3}$ is a vector field on the manifold $SO(3)$.

In order to transform \eqref{secunda} into rotations we need the following computation:
\begin{align*}
\frac{\text{sgn}(\langle {\bf q},{\bf q}_i\rangle)\arccos(|\langle\mathbf{q},\mathbf{q}_{i}\rangle|)}{\sqrt{1-\langle\mathbf{q},\mathbf{q}_{i}\rangle^{2}}} & =\frac{\arccos(|\langle\mathbf{q},\mathbf{q}_{i}\rangle|)\langle {\bf q},{\bf q}_i\rangle}{|\langle {\bf q},{\bf q}_i\rangle|\sqrt{1-\langle\mathbf{q},\mathbf{q}_{i}\rangle^{2}}}=\frac{\frac{|\theta_i|}{2}}{|\cos(\frac{\theta_i}{2})|\sqrt{1-\cos^2(\frac{\theta_i}{2})}}\langle {\bf q},{\bf q}_i\rangle \\
& = \frac{|\theta_i|}{|\sin(\theta_i)|}\langle\mathbf{q},\mathbf{q}_{i}\rangle= \frac{\theta_i}{\sin(\theta_i)}\langle\mathbf{q},\mathbf{q}_{i}\rangle,
\end{align*}
where we have used the property $ {\bf q}\in S^3\backslash \bigcup_{i=1}^r\Pi_i$  which implies that $\langle\mathbf{q},\mathbf{q}_{i}\rangle\neq 0$.
Using \eqref{relatie-delta} we obtain the equivalent system in rotations, and this corresponds to the step {\bf (iv')} of the embedding algorithm,
\begin{equation}
\displaystyle\sum_{i=1}^{r}({\mathbf{R}^{\mathbf{q}_{i}}}^{T}\mathbf{R}^{\mathbf{q}}-
{\mathbf{R}^{\mathbf{q}}}^{T}\mathbf{R}^{\mathbf{q}_{i}})\frac{\theta_{i}}{\sin \theta_{i}}= {\bf 0} \\\Leftrightarrow\sum_{i=1}^{r}\Log({\mathbf{R}^{\mathbf{q}_{i}}}^{T}\mathbf{R}^{\mathbf{q}})={\bf 0}.
\end{equation}
The critical points of the cost function $G_{2_{SO(3)}}$ are the solutions of the equation
\begin{equation}\label{logarithm-equation}
\sum_{i=1}^{r}\Log({\mathbf{R}_i}^{T}\mathbf{R})={\bf 0},
\end{equation}
where $\mathbf{R}\neq \mathbf{R}_i$ and $\theta_i\neq \pm \pi$.
The above equation has been also obtained in \cite{moakher} as a characterization of the
Riemannian mean. This was to be expected, as we are averaging the same cost function by two different methods
that naturally lead to the same result. Also note that the cost function $G_{2_{SO(3)}}$ is not well defined for angles $\theta_i=\pm\pi$ or equivalently $\langle\mathbf{q},\mathbf{q}_{i}\rangle =0$, a situation in which the prolongation function $G_2$ is not differentiable.

Local/global extrema of the function $G_{2_{SO(3)}}:SO(3)\backslash \mathcal{A}\rightarrow \R$ are to be found among the solutions of equation \eqref{logarithm-equation} and the set $\{{\bf R}_1,...,{\bf R}_r\}$. Elements of the set $\{{\bf R}_1,...,{\bf R}_r\}$ can be local/global extrema but this case has to be checked separately because the algorithm breaks down due to the differentiability problems as have been explained above.
\medskip

{\bf III.} Another distance used on the Lie group $SO(3)$ is given by $d_3({\bf R}_1,{\bf R}_2)=1-\frac{1}{2} \sqrt{\text{tr}({\bf R}_1^T{\bf R}_2)+1}$. This distance appears in \cite{huynh}, where the functional dependence between $d_3$ and $d_2$ is also given. We notice that this distance function on $SO(3)$ can be obtained from the pseudodistance $d_3:S^3\times S^3\rightarrow \R$, $d_3({\bf q}_1,{\bf q}_2)=1-|\langle {\bf q}_1,{\bf q}_2\rangle|$. The cost function $G_{3_{SO(3)}}({\bf R})=\sum_{i=1}^r d_3^2({\bf R},{\bf R}_i)$ is lifted to the cost function
$$G_{{3}_{S^{3}}}(\mathbf{q})=\sum_{i=1}^{r}
(1-|\langle\mathbf{q},\mathbf{q}_{i}\rangle|)^{2}.$$

The prolongation function is $G_{3}:\mathbb{R}^{4}\rightarrow\mathbb{R}$, $G_{3}(\mathbf{q})=\sum_{i=1}^{r}
(1-|\langle\mathbf{q},\mathbf{q}_{i}\rangle|)^{2}$. The coefficients of the differential 1-form $dG_3$ are given by:
$$\overline{d G_{3}}(\mathbf{q}) =-2\sum_{i=1}^{r}(1-|\langle\mathbf{q},\mathbf{q}_{i}\rangle|)\text{sgn} (\langle\mathbf{q},\mathbf{q}_{i}\rangle){\bf q}_i,$$
where, in order to have differentiability, we restrict our selves to the open set ${\bf q}\in S^3\backslash \bigcup_{i=1}^r
\Pi_i$ (the set $\Pi_i$ is the hyperplane determined by ${\bf q}_i$ as before).

The system of equations \eqref{sistem} becomes
\begin{equation}
\left\{
  \begin{array}{ll}
    \displaystyle{\sum_{i=1}^{r}}(1-|\langle\mathbf{q},\mathbf{q}_{i}\rangle|)\text{sgn}(\langle\mathbf{q},\mathbf{q}_{i}\rangle)(\langle\mathbf{q},\mathbf{q}\rangle\mathbf{q}_{i}-
\langle\mathbf{q},\mathbf{q}_{i}\rangle\mathbf{q})={\bf 0} \\
    {\bf q}\in S^3\backslash \displaystyle{\bigcup_{i=1}^r}
\Pi_i
  \end{array}
\right.
\end{equation}
and this represents the equation $\mathbf{v_{0}}={\bf 0}$ restricted to the open set $S^3\backslash \bigcup_{i=1}^r
\Pi_i$ of the sphere $S^3$.

As in the previous cases, after applying the projection operator ${\bf P}_*$, we generate the vector field  $\mathbf{P_*}\mathbf{v_{0}}_{|S^3}$ and consequently, we obtain the equivalent system of equations:
\begin{equation}\label{prima}
\left\{
\begin{array}{ll}
    \displaystyle\sum_{i=1}^{r}(1-|\langle\mathbf{q},\mathbf{q}_{i}\rangle|)\text{sgn}(\langle\mathbf{q},\mathbf{q}_{i}\rangle)
\Delta_{i}(\mathbf{q})= {\bf 0} \\
    {\bf q}\in S^3\backslash \displaystyle{\bigcup_{i=1}^r} \Pi_i
\end{array}
\right..
\end{equation}

Transforming the above system into rotations, after some algebraic manipulations using \eqref{relatie-delta} and \eqref{modul}, we obtain the equation from the step {\bf (iv')} of the embedding algorithm
\begin{equation}\label{equation-III}
\sum_{i=1}^r\left ( \frac{2}{\sqrt{\text{tr}({{\bf R}}^T{\bf R}_i)+1}}-1\right ) ({{\bf R}_i}^T{\bf R}-{\bf R}^T{\bf R}_i)={\bf 0},
\end{equation}
where $\theta_i\neq \pm \pi$ or equivalently, ${\bf R}\notin \mathcal{A}=\bigcup_{i=1}^r\mathcal{A}_i$ with
$\mathcal{A}_i=\{{\bf R}\in SO(3)\,|\,\text{tr}({\bf R}_i^T{\bf R})=-1\}.$

Local/global extrema of the function $G_{3_{SO(3)}}:SO(3)\rightarrow \R$ are to be found among the solutions of equation \eqref{equation-III} and the set $\mathcal{A}$. Elements of the set $\mathcal{A}$ can be local/global extrema but this case has to be checked separately because the algorithm breaks down due to the differentiability problems as have been explained above.

\medskip

{\bf  IV.} Next we will apply the techniques developed so far for a cost function that is of $L^p$-type with $p\geq 1$, see \cite{hartly}. We start with the distance $d_1:SO(3)\times
SO(3)\rightarrow\mathbb{R}_{+}$, $d_1(\mathbf{R}_1,\mathbf{R}_2)=||\mathbf{R}_1-\mathbf{R}_2||_{F}$. The associated $L^p$-mean cost function is $G_{(p)_{_{SO(3)}}}:SO(3)\rightarrow \R$,
$$G_{(p)_{_{SO(3)}}}({\bf R})=\sum_{i=1}^r||{\bf R}-{\bf R}_i||_F^p.$$

The lifted cost function is given by
$G_{{(p)}_{S^{3}}}:S^{3}\rightarrow\mathbb{R}$,
$$G_{{(p)}_{S^{3}}}(\mathbf{q})=8^{\frac{p}{2}}\sum_{i=1}^{r}
(1-\langle\mathbf{q},\mathbf{q}_{i}\rangle^{2})^{\frac{p}{2}}.$$

We choose the prolongation function
$G_{(p)}:\mathbb{R}^{4}\rightarrow\mathbb{R}$, $G_{(p)}(\mathbf{q})=8^{\frac{p}{2}}\sum_{i=1}^{r}
(||\mathbf{q}||^2\cdot ||\mathbf{q}_{i}||^2-\langle\mathbf{q},\mathbf{q}_{i}\rangle^{2})^{\frac{p}{2}}$.
By direct computation we obtain:
$$\overline{d G_{(p)}}(\mathbf{q}) =-p8^{\frac{p}{2}}\sum_{i=1}^{r}(||\mathbf{q}||^2\cdot ||\mathbf{q}_{i}||^2-\langle\mathbf{q},\mathbf{q}_{i}\rangle^{2})^{\frac{p}{2}-1} (\langle\mathbf{q},\mathbf{q}_{i}\rangle{\bf q}_i- ||\mathbf{q}_{i}||^2\mathbf{q}).$$
In order to have differentiability of the function $G_{(p)}$, for $p\in [1,2)$ we have to eliminate from the domain of definition the lines $d_i$ that go through points ${\bf 0}$ and $\mathbf{q}_i$.
As in the case of the function $G_1$, the system \eqref{sistem} corresponding to the function $G_{(p)}$ becomes:
\begin{equation}\label{sistem-v0-II}
\left\{\begin{array}{ll}
\displaystyle\sum_{i=1}^{r}\langle\mathbf{q},\mathbf{q}_{i}\rangle(||\mathbf{q}||^2\cdot ||\mathbf{q}_{i}||^2-\langle\mathbf{q},\mathbf{q}_{i}\rangle^2)^{\frac{p}{2}-1}(\langle\mathbf{q},\mathbf{q}\rangle\mathbf{q}_{i}-\langle\mathbf{q},\mathbf{q}_{i}\rangle\mathbf{q})
={\bf 0}
\\<\mathbf{q},\mathbf{q}>=1\\
\mathbf{q}\neq \pm \mathbf{q}_i, \,\,\text{for the case when}\,p\in [1,2)
\end{array}
\right..
\end{equation}

By using the same arguments as before, in the case of function $G_{(p)}$, the equation
$\mathbf{v_{0}}_{|S^{3}}(\mathbf{q})={\bf 0}$ is equivalent with $\mathbf{P_*}\mathbf{v_{0}}_{|S^{3}}={\bf 0}$, which has the following form
\begin{equation}\label{coborata}
\left\{
\begin{array}{ll}
    \displaystyle\sum_{i=1}^{r}(1-\langle\mathbf{q},\mathbf{q}_{i}\rangle^2)^{\frac{p}{2}-1}\langle\mathbf{q},\mathbf{q}_{i}\rangle\Delta_{i}(\mathbf{q})=
    {\bf 0} \\
    \langle\mathbf{q},\mathbf{q}\rangle=1\\
\mathbf{q}\neq \pm \mathbf{q}_i, \,\,\text{for the case when}\,p\in [1,2)
\end{array}
\right..
\end{equation}
Using again \eqref{modul} and \eqref{relatie-delta} we obtain, according to step {\bf (iv')} of the embedding algorithm, the equation in rotations that gives the critical points of the cost function $G_{(p)_{_{SO(3)}}}$,
\begin{equation}\label{equation-p}
\sum_{i=1}^{r}(3-\text{tr}({{\bf R}}^T{\bf R}_i))^{\frac{p}{2}-1}({\mathbf{R}_{i}}^{T}\mathbf{R}-
{\mathbf{R}}^{T}\mathbf{R}_i)={\bf 0},
\end{equation}
with extra condition $\mathbf{R}\neq \mathbf{R}_i$ for the case when $p\in [1,2)$.

Local/global extrema of the function $G_{(p)_{_{SO(3)}}}:SO(3)\rightarrow \R$ are to be found among the solutions of equation \eqref{equation-p} for $p\geq 2$. For the case $p\in [1,2)$ elements of the set $\{{\bf R}_1,...,{\bf R}_r\}$ can be also local/global extrema but this case has to be checked separately because the algorithm breaks down due to the differentiability problems as have been explained above.

\section{Example}

For the $L^p$-type cost function we will analyze the changes that appears in the averaging problem for the case of three sample rotations around $x$-axis when $p=2$, respectively $p=4$.

We consider as sample rotations the following rotations around the $x$-axis:
\begin{equation*}
\mathbf{R}_1=\left(
\begin{array}{ccc}
1 & 0 & 0 \\
0  & -1 & 0 \\
0 & 0 & -1
\end{array}
\right),\,\,
\mathbf{R}_2=\left(
\begin{array}{ccc}
1 & 0 & 0 \\
0  & 0 & -1 \\
0 & 1 & 0
\end{array}
\right),\,\,
\mathbf{R}_3=\left(
\begin{array}{ccc}
1 & 0 & 0 \\
0  & \cos(\alpha) & -\sin(\alpha) \\
0 & \sin(\alpha) & \cos (\alpha)
\end{array}
\right),\,\,\alpha\in [-\pi,\pi],
\end{equation*}
with the corresponding angles of rotation $\theta_1=\pi$, $\theta_2=\frac{\pi}{2}$, and $\theta_3(\alpha)=\alpha$.
For the quaternions associated to this rotations we choose the following:
$${\bf q}_1=(0,1,0,0),\,\,{\bf q}_2=(\frac{\sqrt{2}}{2},\frac{\sqrt{2}}{2},0,0),\,\,{\bf q}_3=(\cos\frac{\alpha}{2},\sin\frac{\alpha}{2},0,0),\,\,\alpha\in [-\pi,\pi].$$

Next, we will find minima/critical points for the cost functions  $G_{(p)_{_{SO(3)}}}$ in the case $p=2$, respectively $p=4$ when the parameter $\alpha$ run between $-\pi$ and $\pi$.  According to the embedding algorithm it is sufficient to study minima/critical points for the corresponding lifted function $G_{{(2)}_{S^{3}}}$, respectively $G_{{(4)}_{S^{3}}}$.
Further, this is equivalent with solving the system of equations \eqref{coborata} for $p=2$ and $p=4$.

We start with the case $p=2$. Solving the system \eqref{coborata} we obtain five families of critical points:
\begin{itemize}
\item [] \hspace{1cm} $\text{Set}_{black}=\{(0,0,\pm\sqrt{1-t^2},t)\,|\,t\in[-1,1]\}$;

\item [] \hspace{1cm} $\text{Set}_{green}=\{(\sqrt{1-x_{2,min}^2(\alpha)},x_{2,min}(\alpha),0,0)\,|\,\alpha\in [-\pi,\pi]\}$;

\item [] \hspace{1cm} $\text{Set}_{pink}=\{(-\sqrt{1-x_{2,min}^2(\alpha)},x_{2,min}(\alpha),0,0)\,|\,\alpha\in [-\pi,\pi]\}$;

\item [] \hspace{1cm} $\text{Set}_{red}=\{(\sqrt{1-x_{2,max}^2(\alpha)},x_{2,max}(\alpha),0,0)\,|\,\alpha\in [-\pi,\pi]\}$;

\item [] \hspace{1cm} $\text{Set}_{blue}=\{(-\sqrt{1-x_{2,max}^2(\alpha)},x_{2,max}(\alpha),0,0)\,|\,\alpha\in [-\pi,\pi]\}$,
\end{itemize}
where $x_{2,min}(\alpha)$ and $x_{2,max}(\alpha)$ are the smallest, respectively largest real positive solutions of the polynomial
\begin{align*}
Q_{2,\alpha}(Z) & =\left( 128 \sin^4 \frac{\alpha}{2}-32
 \sin^2 \frac{\alpha}{2}+4 \right) Z^{4}-\left( 128 \sin^4 \frac{\alpha}{2}-32
 \sin^2 \frac{\alpha}{2}+4 \right) Z^{2} \\
& -16\sin^6 \frac{\alpha}{2}+16\sin^5 \frac{\alpha}{2}\cos \frac{\alpha}{2}+28\sin^4 \frac{\alpha}{2}-8 \sin^2 \frac{\alpha}{2} +1.
\end{align*}
Because the symmetry of the polynomial $Q_{2,\alpha}$ it has only two real positive roots.
\medskip

\begin{figure}[h]
\begin{center}
  \includegraphics[width=15cm]{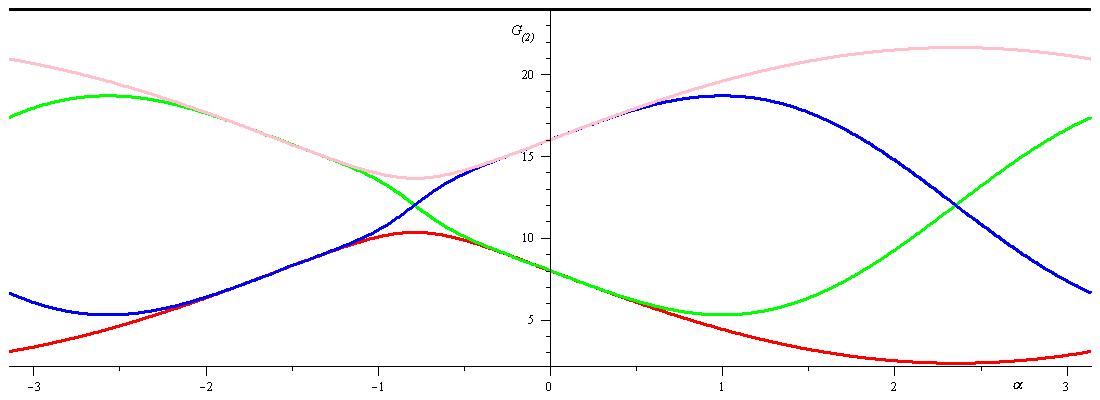}\\
  \caption{The values of  $G_{(2)_{_{SO(3)}}}$ on the critical sets.}\label{fig-G2}
\end{center}
\end{figure}
\medskip

\begin{figure}[h]
\begin{center}
  \includegraphics[width=15cm]{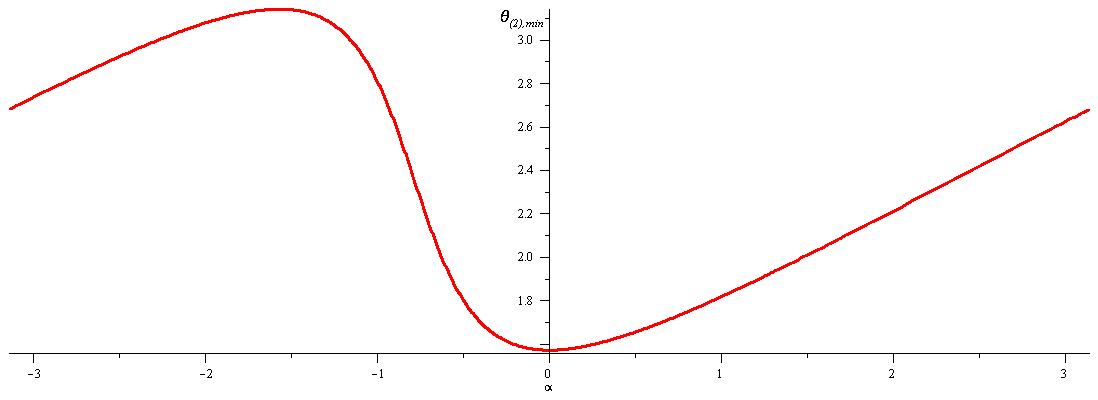}\\
  \caption{The angle of rotation $\theta_{(2),min}$ for rotation where the minimum value is attained.}\label{fig-theta2}
\end{center}
\end{figure}

In the Figure \ref{fig-G2} we represent the values of the cost function $G_{(2)_{_{SO(3)}}}$ on the sets of critical points discovered above when the parameter $\alpha\in [-\pi,\pi]$.  We observe that the absolute maximum is attained on the $\text{Set}_{black}$ which represent rotations of angle $\pi$ around axis that are perpendicular on the $x$-axis and this absolute maximum value does not depend on the parameter $\alpha$. For a fixed value of the parameter $\alpha\in [-\pi,\pi]$ we obtain an absolute minimum on the set $\text{Set}_{red}$ which represent rotations around the $x$-axis. We denote the angle of rotation along the $x$-axis where the minimum value of the cost function $G_{(2)_{_{SO(3)}}}$  is attained  by $\theta_{(2),min}$ and its dependence on the parameter $\alpha$ is represented in the Figure \ref{fig-theta2} bellow.
\medskip

\medskip

Next, we will analyze the behavior of the rotation angle $\theta_{(2),min}$ for different configurations of the sample rotations. For $\alpha=-\pi$ (the rotations ${\bf R}_1$ and ${\bf R}_3$ coincide) we have that the minimum value of the cost function  $G_{(2)_{_{SO(3)}}}$ is attained for a rotation along the $x$-axis with angle $\theta_{(2),min}(-\pi)=2.677945044$. When $\alpha\in [-\pi,-\frac{\pi}{2}]$ we have that $\theta_{(2),min}$ is a monotone increasing function with $\theta_{(2),min}(-\frac{\pi}{2})=\pi$ (minimum rotation coincide with the sample rotation ${\bf R}_1$). For $\alpha\in [-\frac{\pi}{2},0]$, the function $\theta_{(2),min}$ is monotone decreasing and $\theta_{(2),min}(0)=\frac{\pi}{2}$ (minimum rotation coincide with the sample rotation ${\bf R}_2$). If $\alpha\in [0,\pi]$, the function $\theta_{(2),min}$ is again monotone increasing with $\theta_{(2),min}(\pi)=2.677945044$. Consequently, the angle $\theta_{(2),min}$ covers the second quadrant, i.e. $\theta_{(2),min}([-\pi,\pi])=[\frac{\pi}{2},\pi]$.

We will analyze the case $p=4$. The system \eqref{coborata} has a large family of critical points. Among them we obtain as before the $\text{Set}_{black}$ where the absolute maximum for the cost function $G_{(4)_{_{SO(3)}}}$ is attained. In the following, we will study the behavior of the cost function $G_{(4)_{_{SO(3)}}}$ only on the sets of critical points that represent rotations along $x$-axis. The set of critical points that represent rotations along the $x$-axis are constructed, as before from the solutions of the polynomial
$$Q_{4,\alpha}(Z) = a_8Z^8+a_6Z^6+a_4Z^4+a_2Z^2+a_0,$$
where
\begin{itemize}
\item [] $a_8=16$;

\item [] $a_6=-128\sin^4\frac{\alpha}{2}+96\sin^2\frac{\alpha}{2}-128\sin^3\frac{\alpha}{2}\cos\frac{\alpha}{2}+64\sin\frac{\alpha}{2}\cos\frac{\alpha}{2}-32$;

\item [] $a_4=192\sin^4\frac{\alpha}{2}-128\sin^2\frac{\alpha}{2}+192\sin^3\frac{\alpha}{2}\cos\frac{\alpha}{2}-80\sin\frac{\alpha}{2}\cos\frac{\alpha}{2}+24$;

\item []  $a_2=32\sin^6\frac{\alpha}{2}-112\sin^4\frac{\alpha}{2}+48\sin^2\frac{\alpha}{2}-80\sin^3\frac{\alpha}{2}\cos\frac{\alpha}{2}+24\sin\frac{\alpha}{2}\cos\frac{\alpha}{2}-8$;

\item [] $a_0=-16\sin^8\frac{\alpha}{2}+16\sin^6\frac{\alpha}{2}+8\sin^3\frac{\alpha}{2}\cos\frac{\alpha}{2}+1$.
\end{itemize}
Giving values to the parameter $\alpha$, we find the real roots of the polynomial $Q_{4,\alpha}$. For $\alpha\in [-\pi,-1.02)\cup (-0.55,\pi]$ we have two positive real roots $x_{4,min}(\alpha)$ and $x_{4,max}(\alpha)$. As before, we find the four sets of critical points $\text{Set}_{green}$, $\text{Set}_{pink}$, $\text{Set}_{red}$, and $\text{Set}_{blue}$. For $\alpha\in [-1.02,-0.55]$ we find four positive real roots $x_{4,min}(\alpha),\,x_{4,*}(\alpha),\,x_{4,**}(\alpha),\,x_{4,max}(\alpha)$. As sets of critical points, for $\alpha $ in this interval, we find $\text{Set}_{green}$, $\text{Set}_{pink}$, $\text{Set}_{red}$, $\text{Set}_{blue}$, and four more critical sets corresponding to $x_{4,*}(\alpha)$, respectively $x_{4,**}(\alpha)$. Namely, we have the supplementary critical sets corresponding to rotations along the $x$-axis,
\begin{itemize}

\item [] \hspace{1cm} $\text{Set}_{yellow}=\{(\sqrt{1-x_{4,*}^2(\alpha)},x_{4,*}(\alpha),0,0)\,|\,\alpha\in [-1.02,-0.55]\}$;

\item [] \hspace{1cm} $\text{Set}_{violet}=\{(-\sqrt{1-x_{4,*}^2(\alpha)},x_{4,*}(\alpha),0,0)\,|\,\alpha\in [-1.02,-0.55]\}$;

\item [] \hspace{1cm} $\text{Set}_{maroon}=\{(\sqrt{1-x_{4,**}^2(\alpha)},x_{4,**}(\alpha),0,0)\,|\,\alpha\in [-1.02,-0.55]\}$;

\item [] \hspace{1cm} $\text{Set}_{gold}=\{(-\sqrt{1-x_{4,**}^2(\alpha)},x_{4,**}(\alpha),0,0)\,|\,\alpha\in [-1.02,-0.55]\}$.
\end{itemize}

 The function $\alpha\rightarrow x_{4,min}(\alpha)$ is a continuous function on the whole interval $[-\pi,\pi]$, and consequently  $G_{(4)_{_{S^3}}}(\text{Set}_{green})$ and $G_{(4)_{_{S^3}}}(\text{Set}_{pink})$ are continuous curves. The discontinuity of the function  $\alpha\rightarrow x_{4,max}(\alpha)$ implies the discontinuities of the curves $G_{(4)_{_{S^3}}}(\text{Set}_{red})$ and $G_{(4)_{_{S^3}}}(\text{Set}_{blue})$.
\medskip

\begin{figure}[h]
\begin{center}
  \includegraphics[width=15cm]{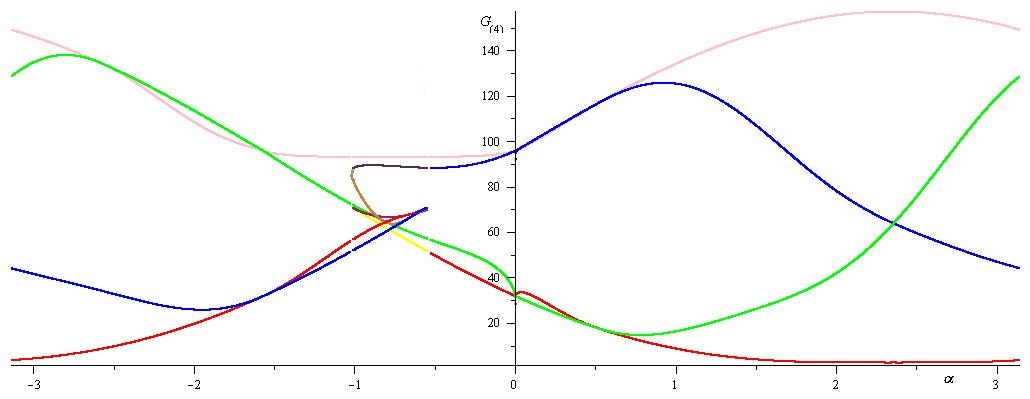}\\
  \caption{The values of  $G_{(4)_{_{SO(3)}}}$ on the critical sets.}\label{fig-G4}
\end{center}
\end{figure}
\medskip

\begin{figure}[h]
\begin{center}
  \includegraphics[width=15cm]{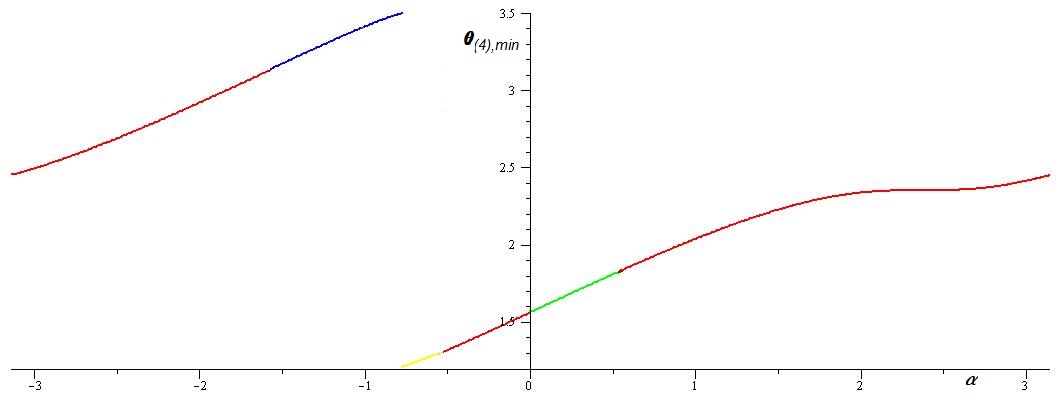}\\
  \caption{The angle  of rotation $\theta_{(4),min}$ for rotations where the minimum value is attained.}\label{fig-theta4}
\end{center}
\end{figure}

In the Figure \ref{fig-G4} we plotted the values of the cost function $G_{(4)_{_{SO(3)}}}$ on the critical sets corresponding to rotations along the $x$-axis when the parameter $\alpha\in [-\pi,\pi]$. In contrast to the case $p=2$, for $p=4$ we have a change of critical sets where the minimum value is attained.  We denote the angle of rotation along the $x$-axis where the minimum value of the function $G_{(4)_{_{SO(3)}}}$  is attained  by $\theta_{(4),min}$ and its dependence on the parameter $\alpha$ is represented in the Figure \ref{fig-theta4} bellow.

\begin{figure}[h]
\begin{center}
  \includegraphics[width=15cm]{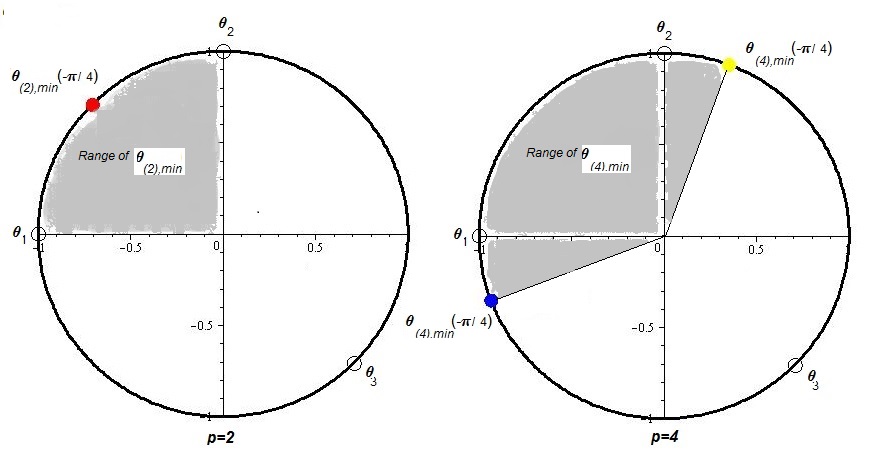}\\
  \caption{Non-uniqueness of the minimal rotation for $\alpha=-\frac{\pi}{4}$ and $p=4$.}\label{cercuri}
\end{center}
\end{figure}

\newpage
Regarding the behavior of the angle $\theta_{(4),min}$ there are several major differences form the behavior of the angle $\theta_{(2),min}$.
\begin{itemize}
\item [(i)] The angle $\theta_{(4),min}$ changes colors, i.e. the minimum value of the cost function $G_{(4)_{_{SO(3)}}}$ is attained for rotations corresponding to quaternions coming from different critical sets having the colors as they appear in the Figure \ref{fig-theta4}.

\item [(ii)] $\theta_{(4),min}$ is a monotone increasing continuous function on each interval $[-\pi,-\frac{\pi}{4})$, respectively $(-\frac{\pi}{4},\pi]$.

\item [(iii)] For $\alpha=-\frac{\pi}{4}$, in contrast to the case $p=2$, the angle $\theta_{(4),min}$ has two values, see Figure \ref{cercuri}, and consequently, $\alpha\rightarrow \theta_{(4),min}(\alpha)$ has two branches. This is due to the non-uniqueness of the critical point where the minimum value of the cost function $G_{(4)_{_{SO(3)}}}$ is attained. More precisely, the critical rotations where the minimum value is attained correspond to the quaternions $(0.82,0.56,0,0)\in \text{Set}_{yellow}$, respectively $(-0.17,0.98,0,0)\in \text{Set}_{blue}$.

\item [(iv)] For the case $p=4$, the angle  $\theta_{(4),min}$ covers the second quadrant and parts of the first and three quadrants, see Figure \ref{fig-theta4} and Figure \ref{cercuri}.

\end{itemize}

\medskip
\newpage

{\bf Acknowledgements.} Petre Birtea and Dan Com\u{a}nescu have been supported by the grant of the Romanian National Authority for
Scientific Research, CNCS UEFISCDI, project number PN-II-RU-TE-2011-3-0006. We are also thankful to Ioan Casu for his help with Maple programming.

\end{document}